\documentclass[12pt]{article}
\usepackage{amssymb,amscd}
\newtheorem{conj}{Conjecture}[section]
\newtheorem{thm}[conj]{Theorem}

\newtheorem{cor}[conj]{Corollary}
\newtheorem{lema}[conj]{Lemma}

\begin{document}
\title{\Large On Perles' question}
\author{{\Large Sini\v sa T. Vre\' cica}
\thanks{Supported by the
Serbian Ministry of Science, Grant 144026.}
\\ {University of Belgrade}}
\maketitle \vskip 2cm

\begin{abstract}
The principal aim of this paper is to determine the minimal
dimension that a nerve of the cover of the sphere $S^h$ by the
open sets not containing a pair of antipodal points could have,
and also to determine the minimal cardinality of such cover (or
the minimal number of vertices of its nerve). In particular, our
result provides the complete answer to the question posed by Micha
Perles (see \cite{st1}). Our results could be seen as the
extensions of the Lyusternik-Schnirel'man version of the
Borsuk-Ulam theorem.

As a consequence, we also obtain the improved lower bound for the
local chromatic number of certain class of graphs.
\end{abstract}

\section{Introduction}

The Lyusternik-Schnirel'man theorem says that if a family of open
sets covers the sphere $S^h$, and none of them contain a pair of
antipodal points, then there are at least $h+2$ of them. It is
easily verified that this statement is one of the many statements
equivalent to the Borsuk-Ulam theorem. A natural question arises
of determining the minimal number of these sets that contain some
point of the sphere. In other words, we consider the question of
determining the minimal dimension that a nerve of such a cover of
the sphere could have. In \cite{st1} and \cite{st2}, G. Simonyi
and G. Tardos gave the complete answer to this question when $h$
is odd, but left the ambiguity of $1$ in the case when $h$ is
even. In this paper, we resolve this ambiguity.

We also discuss the minimal cardinality of the cover of the sphere
satisfying the above properties, or in other words we consider the
minimal number of vertices of the nerve of such a cover of the
sphere of the minimal dimension. In the section 4 we give the
complete answer to this question as well.

The chromatic number of a graph is a very important invariant and
it attracted a lot of attention. Let us mention the nspiring
result of L. Lov\' asz solving the Kneser conjecture and
determining the chromatic number of the certain class of graphs in
terms of connectivity of associated simplicial complex. This
result motivated a lot of research (see \cite{l}).

The local chromatic number of a graph $G$ is defined to be the
minimum number of colors that must appear within distance $1$ of a
vertex of $G$. More formally, if we denote with $N(v)$ the
neighborhood of a vertex $v$ in a graph $G$ (the set of vertices
$v$ is connected to), we define $\psi (G)=\min_c \max_{v\in V(G)}
\vert \{c(u) \mid u\in N(v)\}\vert +1$, where the minimum is taken
over all proper colorings $c$ of $G$, compare \cite{st1}. It is
obvious that the local chromatic number of a graph $G$ is smaller
than its chromatic number. G. Simonyi and G. Tardos showed in
their paper that if a graph $G$ is strongly topologically
$t$-chromatic (which is a little bit stronger assumption than to
be $t$-chromatic), then $\psi (G)\geq \left[\frac t2\right]+1$,
(see \cite{st1}). We don't want to go too much into details here,
and so refer the reader to \cite{st1} for all the definitions,
statements and proofs. As a consequence of our result, we improve
this lower bound for the local chromatic number of graphs by $1$.
There are examples described in \cite{st1} showing that this lower
bound is the best possible.

\section{The minimal intersections of the cover}

Motivated by the related question of Matatyahu Rubin, Micha Perles
formulated the following question.
\bigskip

\noindent {\bf Perles' question.} For which $h$ and $l$, the
sphere $S^h$ can be covered by open sets not containing a pair of
antipodal points in such a way that no point of the sphere is
contained in more than $l$ of these sets?
\bigskip

For given $h$ and $l$, let us denote with $Q(h,l)$ the statement
that the answer to Perles' question is positive. Of course, this
question is closely related to the Lyusternik-Schnirel'man version
of the Borsuk-Ulam theorem. This question could be reformulated so
as to determine the minimal number $l$ such that the sphere $S^h$
can be covered by open sets not containing a pair of antipodal
points such that the intersection of any $l+1$ of them is empty.
This minimal number will be denoted by $Q(h)$.

In two of their papers G. Simonyi and G. Tardos (see \cite{st1}
and \cite{st2}) treated this question. They arrived at this
problem in an attempt to determine the local chromatic number of
certain class of graphs. They were able to prove the following
equivalence.

\begin{thm}
\label{equi} {\bf (\cite{st2})} For every $h$ and $l$ the
statement $Q(h,l)$ is true if and only if there is a continuous
map $g : S^h\rightarrow \| K\|$ to some finite simplicial complex
$K$ of dimension at most $l-1$ satisfying $g(x)\neq g(-x)$ for all
$x\in S^h$. \hfill $\blacksquare$
\end{thm}

They proved a little bit more, namely that one could require the
minimal simplices containing $g(x)$ and $g(-x)$ to be disjoint for
every $x\in S^h$. This was the starting point in their proof that
the statement $Q(h,l)$ is not true if $h\geq 2l-1$. They also
proved that $Q(h,l)$ is true if $h\leq 2l-3$. So, only the case
$h=2l-2$ remained open.

The statement $Q(0,1)$ (obtained when $l=1$) is trivially true,
and they proved (and independently Imre B\' ar\' any) that
$Q(2,2)$ (obtained when $l=2$) is not true.

In terms of the function $Q(h)$, they proved $\frac h2+1\leq
Q(h)\leq \frac h2+2$. So, for odd $h$ we have $Q(h)=\left[\frac
h2\right]+2$ and for even $h$ we have the ambiguity in all cases
except for $h=2$ when we know $Q(2)=3$.

One of the aims of this note is to verify the missing case
$h=2l-2$ completely, i.e. to determine the value $Q(h)$ for all
even $h$. We show the following:

\begin{thm}
\label{main} For every $l\geq 2$ the statement $Q(2l-2,l)$ is not
true, i.e. for even $h$ we have $Q(h)=\left[\frac h2\right]+2$.
\end{thm}

Using the results from \cite{st2}, this theorem implies the
following, improved lower bound for the local chromatic number
$\psi (G)$ of the graphs which are strongly topologically
$t$-chromatic. For the relevant definitions see \cite{st2}.

\begin{cor}
If a graph $G$ is strongly topologically $t$-chromatic, then $\psi
(G)\geq \left[\frac t2\right]+2$. \hfill $\blacksquare$
\end{cor}

The equivalence from the theorem \ref{equi} reduces the proof of
the theorem \ref{main} to the following statement.

\begin{lema}
\label{main2} For every map $g : S^{2k}\rightarrow \| K\|$ to some
finite simplicial complex $K$ of dimension $k$, there is $x\in
S^{2k}$ so that $g(-x)=g(x)$.
\end{lema}

During the preparation of this manuscript we learned that this
lemma is already proved by E.V. \v S\v cepin in \cite{s}. His
proof is based on the fact that the $(2k)$-th homotopy group
$\pi_{2k}\left(\vee_{i=1}^m S_i^k\right)$ of the finite wedge of
$k$-dimensional spheres is finite, and some clever geometrical
constructions. Our proof uses completely different methods, namely
the Smith theory and cohomological methods. We use the similar
methods in the final section of the paper, and so we present our
proof of the lemma \ref{main2} in the next section.

\section{Proof of the main lemma}

Suppose, to the contrary, that for some map $g$ such a point $x\in
S^{2k}$ does not exist. The map $g$ induces the map $G :
S^{2k}\rightarrow \| K\| \times \| K\|$ defined by
$G(x)=(g(x),g(-x))$, which misses the diagonal $\Delta =\{(u,u)
\mid u\in \| K\|\}\subseteq \| K\| \times \| K\|$. Let us denote
with $X=(\| K\| \times \| K\|)\setminus \Delta$. So, we have the
map $F : S^{2k}\rightarrow X$ defined by $F(x)=(g(x),g(-x))$. The
group $\mathbb{Z}/2$ acts naturally and freely on these spaces and
the map $F$ is easily seen to be equivariant.

Throughout this paper all chain complexes and homology and
cohomology groups are considered with $\mathbb{Z}/2$ coefficients
which are dropped from the notation. For each free
$(\mathbb{Z}/2)$-space $A$, there is the Smith exact sequence of
chain complexes and chain maps:
$$0\longrightarrow C^s_{\ast}(A)\stackrel{i}{\longrightarrow} C_{\ast}(A)
\stackrel{q}\longrightarrow C^a_{\ast}(A)\longrightarrow 0,$$

\noindent where the symmetric chains $C^s_{\ast}(A)$ are those
chains $c$ satisfying $\sigma c=c$ for the generator $\sigma$ of
the group $\mathbb{Z}/2$, and the antisymmetric chains
$C^a_{\ast}(A)$ are the chains satisfying $\sigma (c)=-c$.
(Consult \cite{b}.) The chain map $i$ is the inclusion of the
symmetric chains in the chain complex of all chains, and the chain
map $q$ is defined by $q(c)=c-\sigma (c)$. Moreover this exact
sequence is natural, and so for the equivariant map $F :
S^{2k}\rightarrow X$ of free $(\mathbb{Z}/2)$-spaces we have a
commutative diagram:

$$
\begin{CD}
0 @>>> C^s_{\ast}(S^{2k}) @>i>> C_{\ast}(S^{2k}) @>q>> C^a_{\ast}(S^{2k}) @>>> 0 \\
@. @VV{F_s}V @VVFV @VV{F_a}V @. \\
0 @>>> C^s_{\ast}(X) @>i>> C_{\ast}(X) @>q>> C^a_{\ast}(X) @>>> 0
\end{CD}
$$

This diagram induces the commutative diagram of the long exact
sequences in cohomology. Since the spaces $S^{2k}$ and $X$ are
$2k$-dimensional, the cohomology groups in the dimension $2k+1$
and higher are trivial and the part of the diagram has the
following form.

$$
\begin{CD}
\cdots @>>> H_a^{2k}(S^{2k}) @>{q^{\ast}}>> H^{2k}(S^{2k}) @>{i^{\ast}}>> H_s^{2k}(S^{2k}) @>>> 0 \\
@. @AA{F_a^{\ast}}A @AA{F^{\ast}}A @AA{F_s^{\ast}}A @. \\
\cdots @>>> H_a^{2k}(X) @>{q^{\ast}}>> H^{2k}(X) @>{i^{\ast}}>>
H_s^{2k}(X) @>>> 0
\end{CD}
$$

Since the action of the group $\mathbb{Z}/2$ on the spaces
$S^{2k}$ and $X$ is free, the cohomology ring of symmetric chain
complexes is isomorphic to the usual cohomology ring of the spaces
of orbits $S^{2k}/(\mathbb{Z}/2)=\mathbb{R}P^{2k}$ and
$X/(\mathbb{Z}/2)$. Moreover, the homomorphism $F_s^{\ast} :
H^1(X/(\mathbb{Z}/2))\rightarrow H^1(\mathbb{R}P^{2k})$ sends the
generator $v\in H^1(X/(\mathbb{Z}/2))$ to the generator $u\in
H^1(\mathbb{R}P^{2k})$. Since $F_s^{\ast}$ is a ring homomorphism,
we have $F_s^{\ast}(v^{2k})=u^{2k}\neq 0$ in
$H^{2k}(\mathbb{R}P^{2k})$.

From the exact sequence in the above diagram, it follows that the
homomorphism $i^{\ast} : H^{2k}(X)\rightarrow
H^{2k}(X/(\mathbb{Z}/2))$ is onto, and so there is $w\in
H^{2k}(X)$ so that $i^{\ast}(w)=v^{2k}$. From the commutativity of
the above diagram we have
$i^{\ast}(F^{\ast}(w))=F_s^{\ast}(i^{\ast}(w))=u^{2k}\neq 0$.

Finally, we reach a contradiction by proving $F^{\ast}(w)=0$. Let
$j : X\rightarrow \|K\| \times \|K\|$ be the inclusion map. From
the long exact sequence in cohomology of the pair $(\|K\| \times
\|K\|,X)$ it follows that the map $j^{\ast} : H^{2k}(\|K\| \times
\|K\|)\rightarrow H^{2k}(X)$ is onto, since $H^{2k+1}(\|K\| \times
\|K\|,X)=0$ for the dimensional reasons. So, it is enough to prove
that the map $F^{\ast}\circ j^{\ast}$ is trivial.

By the K\" unneth formula and since the cross product could be
expressed by $a\times b=p_1^{\ast}(a)\cup p_2^{\ast}(b)$ (where
$p_1,p_2 : \|K\| \times \|K\|\rightarrow \|K\|$ are projections),
any element in $H^{2k}(\|K\| \times \|K\|)$ could be expressed as
the cup product of two elements from $H^k(\|K\| \times \|K\|)$.
(Consult \cite{dk} or \cite{h}. Note that $\dim K=k$ and that we
work with $\mathbb{Z}/2$ coefficients.) Since $F^{\ast}\circ
j^{\ast}$ is a ring homomorphism, the image of any element from
$H^{2k}(\|K\| \times \|K\|)$ is the cup product of two elements
from $H^k(S^{2k})$. But the latter group is trivial and so is the
map $F^{\ast}\circ j^{\ast}$. The proof follows. \hfill
$\blacksquare$
\bigskip

\section{The minimal cardinality of the cover}

In this section we address the question of the determination of
the minimal total number of sets $n$ of some cover of the sphere
with the above properties, or in other words we consider the
question of the determination of the minimal number of vertices of
the nerve of this cover of minimal dimension. It is obvious that
the sphere $S^1$ has a covering with three sets so that any point
in the sphere is contained in at most two sets. Also, the sphere
$S^2$ has a covering with four sets so that any point in the
sphere is contained in at most three sets from the covering. For
the higher dimensional spheres the question is more complicated,
but the answers are more regular. We also discuss the reason for
the answer in the cases of the one- and two-dimensional sphere to
be different from the general scheme. We start by proving a
geometrical lemma.

\begin{lema}
If the sphere $S^h$ could be covered with $n$ sets which do not
contain a pair of antipodal points such that every point of the
sphere is contained in at most $l$ sets, then the sphere $S^{h+1}$
could be covered by $n+1$ sets which do not contain a pair of
antipodal points such that every point of the sphere is contained
in at most $l+1$ sets.
\end{lema}

\medskip\noindent
{\bf Proof:} Let us consider the sphere $S^h$ as the equatorial
sphere of the sphere $S^{h+1}$, and its covering
$\mathcal{U}=\{U_1,...,U_n\}$ with the required properties. We
describe the covering
$\mathcal{V}=\{V_1,...,V_{n-1},V_n,V_{n+1}\}$ of the sphere
$S^{h+1}$.

Let $\varepsilon$ be some small positive number. For
$i=1,2,...,n-1$, let $V_i$ be the union of all arcs on the great
circle containing the northern and the southern pole, and
containing the points of $U_i$, and going from each of these
points $4\varepsilon$ in the direction of the northern pole, and
$2\varepsilon$ in the direction of the southern pole. Let $V_n$ be
the union of the set formed in the same way from the set $U_n$,
and the upper hemisphere without the $3\varepsilon$ thick closed
neighborhood of the equatorial sphere. Finally, let $V_{n+1}$ be
the lower hemisphere without the $\varepsilon$ thick closed
neighborhood of the equatorial sphere.

The reader will easily verify that the covering $\mathcal{V}$ of
the sphere $S^{h+1}$ satisfies the required properties. \hfill
$\blacksquare$
\bigskip

In the next two theorems we determine the minimal number of sets
in the cover of the sphere with the required properties. The first
theorem claims that a cover of given cardinality exists, and has
also a geometrical proof. The second theorem states that the given
cardinality is minimal and its proof is topological.

\begin{thm}
The odd-dimensional sphere $S^{2k-1}$ could be covered with $2k+2$
sets which do not contain a pair of antipodal points such that
every point of the sphere is contained in at most $k+1$ sets of
the cover. The even-dimensional sphere $S^{2k}$ could be covered
with $2k+3$ sets which do not contain a pair of antipodal points
such that every point of the sphere is contained in at most $k+2$
sets of the cover.
\end{thm}

\medskip\noindent
{\bf Proof:} The first statement of the theorem is proved in
\cite{st1}, and the second statement follows from the first one
and from the preceding lemma. \hfill $\blacksquare$
\bigskip

Let us denote with $K_{\delta}^2$ the deleted square (here we
consider the deleted product of simplicial complexes) of the
complex $K$ (i.e. the regular CW-complex whose cells are the
products of two disjoint faces of $K$).

\begin{thm}
For the sphere of the dimension at least $3$ the estimate obtained
in the previous theorem is the best possible.
\end{thm}

\medskip\noindent
{\bf Proof:} Let us first prove that the odd dimensional sphere
$S^{2k-1}$ (when $k\geq 2$) could not be covered with $2k+1$ sets
not containing a pair of antipodal points, such that every point
of the sphere is contained in at most $k+1$ sets of the cover.
According to the theorem \ref{equi}, it suffices to prove that
there is no mapping from the sphere $S^{2k-1}$ to the
$k$-dimensional skeleton $K$ of the $(2k)$-dimensional simplex
$\sigma =(v_0,v_1,...,v_{2k})$, mapping every pair of the
antipodal points of the sphere into disjoint simplices of $K$.
(Note also the remark after the theorem \ref{equi}.)

If such a map would exist, it would induce an
$(\mathbb{Z}/2)$-equivariant map $\varphi : S^{2k-1}\rightarrow
K_{\delta}^2$ to the deleted square of the complex $K$. The map
$\varphi$ would induce the map of the spaces of orbits
$\tilde{\varphi} : \mathbb{R}P^{2k-1}\rightarrow
K_{\delta}^2/(\mathbb{Z}/2)$. The dimension of these complexes is
$2k-1$, and their top-dimensional homology groups are trivial.
Namely, every top-dimensional cell (e.g. $(v_0,v_1,...,v_k)\times
(v_{k+1},...,v_{2k-1},v_{2k})$) has a facet (e.g.
$(v_0,v_1,...,v_k)\times (v_{k+1},...,v_{2k-1})$) which is
contained only in this top-dimensional cell. The same is true for
the regular CW-complex of the space of orbits, and so
$H_{2k-1}(K_{\delta}^2/(\mathbb{Z}/2))=0$. From this, it follows
$H^{2k-1}(K_{\delta}^2/(\mathbb{Z}/2))=0$.

Now it is not difficult to obtain a contradiction in a similar way
as in the preceding section. Namely, the mapping $\tilde
{\varphi}^{\ast} : H^1(K_{\delta}^2/(\mathbb{Z}/2))\rightarrow
H^1(\mathbb{R}P^{2k-1})$ is non-trivial, i.e. it maps a generator
$v$ to the generator $u$. So, $\tilde
{\varphi}^{\ast}(v^{2k-1})=u^{2k-1}\neq 0$. But, $v^{2k-1}\in
H^{2k-1}(K_{\delta}^2)=0$, and we reached a contradiction.

Now we prove that the even dimensional sphere $S^{2k}$ (when
$k\geq 2$) could not be covered with $2k+2$ sets not containing a
pair of antipodal points, such that every point of the sphere is
contained in at most $k+2$ sets of the cover. Again, according to
the theorem \ref{equi}, it suffices to prove that there is no
mapping from the sphere $S^{2k}$ to the $(k+1)$-dimensional
skeleton $K$ of the $(2k+1)$-dimensional simplex $\sigma
=(v_0,v_1,...,v_{2k},v_{2k+1})$ mapping every pair of the
antipodal points of the sphere into disjoint simplices of $K$.

If such a map would exist, it would induce an
$(\mathbb{Z}/2)$-equivariant map $\varphi : S^{2k}\rightarrow
K_{\delta}^2$. The map $\varphi$ would induce the map of the
spaces of orbits $\tilde{\varphi} : \mathbb{R}P^{2k}\rightarrow
K_{\delta}^2/(\mathbb{Z}/2)$. The dimension of these complexes is
$2k$, and their top-dimensional homology groups are trivial.
Namely, there are two types of the top-dimensional cells in
$K_{\delta}^2$ and its orbit space $K_{\delta}^2/(\mathbb{Z}/2)$,
and those are the cells of the form $(v_0,v_1,...,v_k)\times
(v_{k+1},...,v_{2k},v_{2k+1})$, and the cells of the form
$(v_0,v_1,...,v_k,v_{k+1})\times (v_{k+2},...,v_{2k},v_{2k+1})$.
All the facets of the cells of the first type are also the facets
of another top-dimensional cell, and this cell is of the second
type. Every cell of the second type has a facet (of the form
$(v_0,v_1,...,v_k,v_{k+1})\times (v_{k+2},...,v_{2k})$) which is
contained only in this top-dimensional cell. The same is true for
the regular CW-complex of the space of orbits, and so
$H_{2k}(K_{\delta}^2/(\mathbb{Z}/2))=0$. From this, it follows
$H^{2k}(K_{\delta}^2/(\mathbb{Z}/2))=0$. Now, we reach a
contradiction in the same way as for odd-dimensional spheres.
\hfill $\blacksquare$
\bigskip

\noindent {\bf Remark 4.4} If $k=1$ in the case of the
odd-dimensional sphere (i.e. if we consider the sphere $S^1$), the
facet contained in only one top-dimensional cell does not exist,
since in this case it would be of type $(v_0,v_1)\times \emptyset$
which is the empty set. Similarly, if $k=1$ in the case of the
even-dimensional sphere (i.e. if we consider the sphere $S^2$),
the facet of the top-dimensional cell of the second type contained
only in this top-dimensional cell would be of type
$(v_0,v_1,v_2)\times \emptyset$ which is the empty set.

So, in the cases of the $1$- and $2$-dimensional spheres, the
argument does not work. Moreover, as we mentioned, it is not
difficult to construct the coverings in these cases with smaller
number of sets. The sets in these coverings are $\varepsilon$
neighborhoods of facets of a simplex inscribed in these spheres.

\vskip 1cm

\parbox{6cm}{Sini\v sa T. Vre\' cica \par Faculty of Mathematics
\par University of Belgrade
\par Studentski trg 16, P.O.B. 550 \par
\par 11000 Belgrade \par vrecica@matf.bg.ac.yu}

\end{document}